# THE BROWNIAN WEB: CHARACTERIZATION AND CONVERGENCE

By L. R. G. Fontes[1], M. Isopi[2], C. M. Newman[3] and
K. Ravishankar[4]

*Universidade de São Paulo, Università di Roma, New York University and
SUNY College at New Paltz*

The Brownian web (BW) is the random network formally consisting of the paths of coalescing one-dimensional Brownian motions starting from every space-time point in $\mathbb{R} \times \mathbb{R}$. We extend the earlier work of Arratia and of Tóth and Werner by providing a new characterization which is then used to obtain convergence results for the BW distribution, including convergence of the system of all coalescing random walks to the BW under diffusive space-time scaling.

**1. Introduction.** In this paper, we present a number of results concerning the characterization of and convergence to a striking stochastic object called the *Brownian web* (BW). Several of the main results were previously announced, with sketches of the proofs, in [13].

Roughly speaking, the BW is the collection of graphs of coalescing one-dimensional Brownian motions (with unit diffusion constant and zero drift) starting from all possible starting points in one plus one-dimensional (continuous) space-time. This object was originally studied more than twenty years ago by Arratia [5], motivated by asymptotics of one-dimensional voter models, and then about five years ago by Tóth and Werner [26], motivated by the problem of constructing continuum "self-repelling motions." Our own interest in this object arose because of its relevance to "aging" in statistical physics models of one-dimensional coarsening [14, 15]—which returns us to

Received April 2003; revised November 2003.
[1]Supported in part by CNPq Grants 300576/92-7 and 662177/96-7 (PRONEX) and FAPESP Grant 99/11962-9.
[2]Supported in part by the CNPq-CNR agreement.
[3]Supported in part by NSF Grant DMS-01-04278.
[4]Supported in part by NSF Grant DMS-98-03267.

Arratia's original context of voter models, or equivalently coalescing random walks in one dimension. This motivation leads to our primary concern with weak convergence results, which in turn requires a careful choice of space for the BW so as to obtain useful characterization criteria for its distribution.

We remark that there are two questions we do not address in this paper that are worthy of consideration. The first is whether our convergence results might play some role in studying the convergence of discrete to continuum self-repelling motion [26]. The second is whether there are interesting connections between the BW and super-Brownian motions; in this regard, the work of [8, 9] may be relevant since it deals with noncrossing paths.

We continue the Introduction by discussing coalescing random walks and their scaling limits. Let us begin by constructing random paths in the plane, as follows. Consider the two-dimensional lattice of all points $(i,j)$ with $i,j$ integers and $i+j$ even. Let a walker at spatial location $i$ at time $j$ move right or left at unit speed between times $j$ and $j+1$ if the outcome of a fair coin toss is heads ($\Delta_{i,j} = +1$) or tails ($\Delta_{i,j} = -1$), with the coin tosses independent for different space-time points $(i,j)$. Figure 1 depicts a simulation of the resulting paths.

The path of a walker starting from $y_0$ at time $s_0$ is the graph of a simple symmetric one-dimensional random walk, $Y_{y_0,s_0}(t)$. At integer times, $Y_{y_0,s_0}(t)$ is the solution of the simple stochastic difference equation,

(1.1) $$Y(j+1) - Y(j) = \Delta_{Y(j),j}, \qquad Y(s_0) = y_0.$$

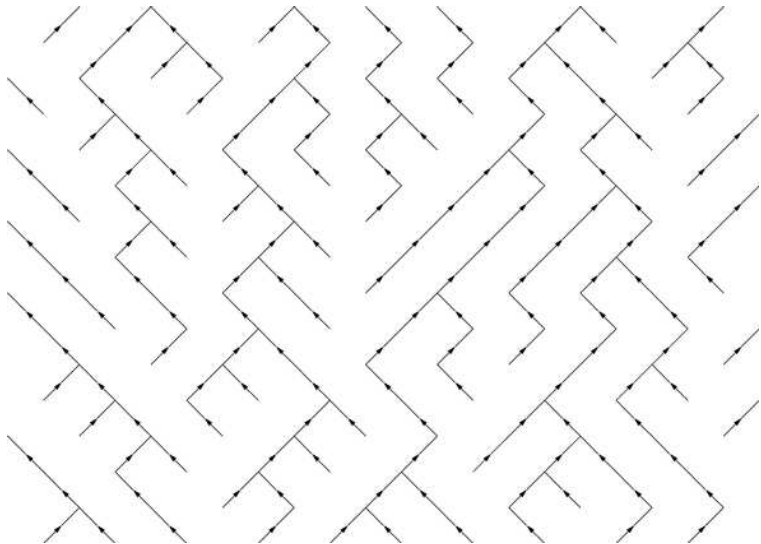

FIG. 1. *Coalescing random walks in discrete time; the horizontal coordinate is space and the vertical coordinate is time.*

Furthermore, the paths of distinct walkers starting from different $(y_0, s_0)$'s are automatically *coalescing*—that is, they are independent of each other until they coalesce (i.e., become identical) upon meeting at some space-time point.

If the increments $\Delta_{i,j}$ remain i.i.d. but take values besides $\pm 1$ (e.g., $\pm 3$), then one obtains nonsimple random walks whose paths can cross each other in space-time, although they still coalesce when they land on the same space-time lattice site. Such systems with crossing paths will be discussed in Section 5 (see also [20]).

After rescaling to spatial steps of size $\delta$ and time steps of size $\delta^2$, a single rescaled random walk (say, starting from 0 at time 0) $Y_{0,0}^{(\delta)}(t) = \delta Y_{0,0}(\delta^{-2} t)$ converges as $\delta \to 0$ to a standard Brownian motion $B(t)$. That is, by the Donsker invariance principle [10], the distribution of $Y_{0,0}^{(\delta)}$ on the space of continuous paths converges as $\delta \to 0$ to standard Wiener measure.

The invariance principle is also valid for continuous-time random walks, where the move from $i$ to $i \pm 1$ takes an exponentially distributed time. In continuous time, coalescing random walks are at the heart of Harris' graphical representation of the (one-dimensional) voter model [18] and their scaling limits arise naturally in the physical context of (one-dimensional) aging (see, e.g., [14, 15]). Of course, finitely many rescaled coalescing walks in discrete or continuous time (with rescaled space-time starting points) converge in distribution to finitely many coalescing Brownian motions. In this paper, we present results concerning the convergence in distribution of the complete collection of the rescaled coalescing walks from *all* the starting points.

Our results are in two main parts:

1. A new characterization of the limiting object, the standard BW.
2. Convergence criteria, which are applied, in this paper, to coalescing random walks.

As a cautionary remark, we point out that the scaling limit motivating our convergence results does not belong to the realm of hydrodynamic limits of particle systems but rather to the realm of invariance principles.

A key ingredient of the new characterization and the convergence is the choice of a space on which the BW measure is defined; this is a space whose elements are *collections* of paths (see Sections 2 and 3). The convergence criteria and application (see, e.g., Theorems 2.2 and 6.1) are the BW analogues of Donsker's invariance principle. Like Brownian motion itself, we expect that the BW and its variants will be quite ubiquitous as scaling limits, well beyond the context of coalescing random walks (and our sufficient conditions for convergence). One situation where this occurs is for two-dimensional "Poisson webs" [11]. Another example is in the area of river basin modelling; in [24] (see also [16, 21, 27]), coalescing random walks were proposed as a model of a drainage network. Some of the questions about scaling in such models may find answers in the context of their scaling limits. For more on coalescing random walk and other models for river basins, see [23, 30, 31].

Much of the construction of the BW (but without convergence results) was already done in the groundbreaking work of Arratia [4, 5] (see also [6, 19]) and then in the work of Tóth and Werner [26] who derived many important properties of the BW (see also [25] and [28]; in the latter reference, the BW is introduced in relation to *black noise*). Arratia, Tóth and Werner all recognized that in the limit $\delta \to 0$ there would be (nondeterministic) space-time points $(x,t)$ starting from which there are multiple limit paths and they provided various conventions (e.g., semicontinuity in $x$) to avoid such multiplicity. An important feature of our approach to the BW is to accept the intrinsic nonuniqueness by choosing an appropriate metric space in which the BW takes its values. Roughly speaking, instead of using some convention to obtain a process that is a *single-valued* mapping from each space-time starting point to a single path from that starting point, we allow *multivalued* mappings; more accurately, our BW value is the collection of *all* paths from all starting points. This choice of space is very much in the spirit of earlier work [1, 2, 3] on spatial scaling limits of critical percolation models and spanning trees, but modified for our particular space-time setting; the directed (in time) nature of our paths considerably simplifies the topological setting compared to [1, 2, 3].

The Donsker invariance principle implies that the distribution of any continuous (in the sup-norm metric) functional of $Y_{0,0}^{(\delta)}$ converges to that for Brownian motion. The classic example of such a functional is the random walk maximum, $\sup_{0 \leq t \leq 1} Y_{0,0}^{(\delta)}(t)$. An analogous example for coalescing random walks is the maximum over all rescaled walks starting at (or passing through) some vertical (time-like) interval, that is, the maximum value (for times $t \in [s,1]$) over walks touching any space-time point of the form $(0,s)$ for some $s \in [0,1]$. In this case, the functional is not quite continuous for our choice of metric space, but it is continuous almost everywhere (with respect to the BW measure), which is sufficient.

The rest of the paper is organized as follows. Section 2 contains two theorems. The first, Theorem 2.1, is a characterization of the BW, as in [5, 26] but adapted to our choice of space; the second and one of our main results is Theorem 2.2 which is a convergence theorem for the important special case where, even before taking a limit, all paths are noncrossing. Section 3 contains propositions related to Theorem 2.1, as well as an alternative characterization, Theorem 3.1, in which a kind of separability condition is replaced by a minimality condition. In Section 4, we present our new characterization results (Theorems 4.1 and 4.2) based on certain counting random variables, which will be needed for the derivation of our main convergence results. We remark that there are analogous characterization and convergence results jointly for the BW and its dual web of backward paths (important properties of the BW and its dual may be found in [25, 26]; see also [12]). In Section 5, we extend our convergence results to cover the case of crossing paths; the proof of the noncrossing result, Theorem 2.2, is given here as a corollary of the more general result. In Section 6, we apply our (noncrossing) convergence results to the case of coalescing random walks. There are two appendices: the first covers issues of measurability, the second issues of compactness and tightness.

A number of theorems and propositions in this paper are either essentially contained in or easily derived from [5] and/or [26]. In those cases, we have omitted the proofs and simply refer the reader to the papers cited above. Detailed proofs can be found in a previous longer version of this paper [12], which uses the same choice of space and notation as this paper.

**2. Convergence for noncrossing paths.** We begin with three metric spaces: $(\overline{\mathbb{R}}^2, \rho)$, $(\Pi, d)$ and $(\mathcal{H}, d_{\mathcal{H}})$. The elements of the three spaces are, respectively: points in space-time, paths with specified starting points in space-time and collections of paths with specified starting points. The BW will be an $(\mathcal{H}, \mathcal{F}_{\mathcal{H}})$-valued random variable, where $\mathcal{F}_{\mathcal{H}}$ is the Borel $\sigma$-field associated to the metric $d_{\mathcal{H}}$. Complete definitions of these three spaces will be given in Section 3.

For an $(\mathcal{H}, \mathcal{F}_{\mathcal{H}})$-valued random variable $\overline{\mathcal{W}}$ (or its distribution $\mu$), we define the *finite-dimensional distributions* of $\overline{\mathcal{W}}$ as the induced probability measures $\mu_{(x_1,t_1;\ldots;x_n,t_n)}$ on the subsets of paths starting from any finite deterministic set of points $(x_1, t_1), \ldots, (x_n, t_n)$ in $\mathbb{R}^2$. There are several ways in which the BW can be characterized; they differ from each other primarily in the type of extra condition required beyond the finite-dimensional distributions (which are those of coalescing Brownian motions). In the next theorem, the extra condition is a type of Doob separability property (see, e.g., Chapter 3 of [29]). Variants are stated later either using a minimality property (Theorem 3.1) or a counting random variable (Theorems 4.1 and 4.2). Theorem 4.2 is the one most directly suited to the convergence results of Section 5.

The events and random variables appearing in the next two theorems are $(\mathcal{H}, \mathcal{F}_{\mathcal{H}})$-measurable. This claim follows straightforwardly from Proposition A.1. The proof of Theorem 2.1 follows primarily from Propositions 3.1 and 3.3. The proofs of these propositions are essentially contained [5] and [26].

THEOREM 2.1. *There is an $(\mathcal{H}, \mathcal{F}_{\mathcal{H}})$-valued random variable $\overline{\mathcal{W}}$ whose distribution is uniquely determined by the following three properties:*

(o) *From any deterministic point $(x,t)$ in $\mathbb{R}^2$, there is almost surely a unique path $W_{x,t}$ starting from $(x,t)$.*

(i) *For any deterministic $n, (x_1, t_1), \ldots, (x_n, t_n)$, the joint distribution of $W_{x_1,t_1}, \ldots, W_{x_n,t_n}$ is that of coalescing Brownian motions (with unit diffusion constant).*

(ii) *For any deterministic, dense countable subset $\mathcal{D}$ of $\mathbb{R}^2$, almost surely, $\overline{\mathcal{W}}$ is the closure in $(\mathcal{H}, d_{\mathcal{H}})$ of $\{W_{x,t} : (x,t) \in \mathcal{D}\}$.*

REMARK 2.1. One can choose a single dense countable $\mathcal{D}_0$ and in (o), (i) and (ii) restrict to space-time starting points from that $\mathcal{D}_0$. Different characterization theorems for the BW with alternatives for (ii) are given in Sections 3 and 4. We note that there are natural $(\mathcal{H}, \mathcal{F}_{\mathcal{H}})$-valued random variables satisfying (o) and (i) but not (ii). An instance of such a random variable will be studied elsewhere, and shown to arise as the scaling limit of stochastic flows, extending earlier work of Piterbarg [22].

The next theorem is our convergence result for noncrossing processes; a more general result is given in Section 5. We first define a counting variable essential to all our convergence results and the related new characterization results of Section 4.

DEFINITION 2.1. For $t > 0$ and $t_0, a, b \in \mathbb{R}$ with $a < b$, let $\eta(t_0, t; a, b)$ be the number of *distinct* points in $\mathbb{R} \times \{t_0 + t\}$ that are touched by paths in $\overline{\mathcal{W}}$ which also touch some point in $[a,b] \times \{t_0\}$. Let also $\hat{\eta}(t_0, t; a, b) = \eta(t_0, t; a, b) - 1$.

We note that by duality arguments (see [5, 26]), it can be shown that for deterministic $t_0, t, a, b$, this $\hat{\eta}$ is equidistributed with the number of distinct points in $[a,b] \times \{t_0 + t\}$ that are touched by paths in $\overline{\mathcal{W}}$ which also touch $\mathbb{R} \times \{t_0\}$.

THEOREM 2.2. *Suppose $\mathcal{X}_1, \mathcal{X}_2, \ldots$ are $(\mathcal{H}, \mathcal{F}_{\mathcal{H}})$-valued random variables with noncrossing paths. If, in addition, the following three conditions are valid, then the distribution $\mu_n$ of $\mathcal{X}_n$ converges to the distribution $\mu_{\overline{W}}$ of the standard BW.*

(I1) *There exist $\theta_n^y \in \mathcal{X}_n$ for $y \in \mathbb{R}^2$ satisfying: for any deterministic $y_1, \ldots, y_m \in \mathcal{D}$, $\theta_n^{y_1}, \ldots, \theta_n^{y_m}$ converge in distribution as $n \to \infty$ to coalescing Brownian motions (with unit diffusion constant) starting at $y_1, \ldots, y_m$.*

(B1) *$\forall t > 0$, $\limsup_{n \to \infty} \sup_{(a,t_0) \in \mathbb{R}^2} \mu_n(\hat{\eta}(t_0, t; a, a+\varepsilon) \geq 1) \to 0$ as $\varepsilon \to 0^+$.*

(B2) *$\forall t > 0, \varepsilon^{-1} \limsup_{n \to \infty} \sup_{(a,t_0) \in \mathbb{R}^2} \mu_n(\hat{\eta}(t_0, t; a, a+\varepsilon) \geq 2) \to 0$ as $\varepsilon \to 0^+$.*

Convergence of coalescing random walks (in discrete and continuous time) (see Theorem 6.1) is obtained as a corollary to Theorem 2.2.

**3. Construction and initial characterizations.** In this section, we discuss the construction of the BW and the proof of Theorem 2.1. Then we give in Theorem 3.1 a somewhat different characterization of the BW distribution.

Let $(\Omega, \mathcal{F}, \mathbb{P})$ be a probability space where an i.i.d. family of standard Brownian motions $(B_j)_{j \geq 1}$ is defined. Let $\mathcal{D} = \{(x_j, t_j), j \geq 1\}$ be a countable dense set in $\mathbb{R}^2$. Let $W_j$ be a Brownian path starting at position $x_j$ at time $t_j$. More precisely,

$$(3.1) \qquad W_j(t) = x_j + B_j(t - t_j), \qquad t \geq t_j.$$

We now construct, following [5], coalescing Brownian paths out of the family of paths $(W_j)_{j \geq 1}$ by specifying coalescing rules. When two paths meet for the first time, they coalesce into a single path, which is that of the Brownian motion with the lower index. We denote the coalescing Brownian paths by $\widetilde{W}_j, j \geq 1$. Notice that the strong Markov property of Brownian motion allows for a lot of freedom in giving a coalescing rule. Any rule, even nonlocal, that does not depend on the realization of the $(W_j)$'s *after* the time of coalescence will yield the same object in distribution. General definitions and constructions of coalescing Brownian motions can be found in [5].

We define the BW skeleton $\mathcal{W}(\mathcal{D})$ with *starting set* $\mathcal{D}$ by

$$(3.2) \qquad \mathcal{W}_k = \mathcal{W}_k(\mathcal{D}) = \{\widetilde{W}_j : 1 \leq j \leq k\},$$

$$(3.3) \qquad \mathcal{W} = \mathcal{W}(\mathcal{D}) = \bigcup_k \mathcal{W}_k.$$

Now we give detailed definitions of the three spaces introduced in Section 2. $(\overline{\mathbb{R}}^2, \rho)$ is the completion (or compactification) of $\mathbb{R}^2$ under the metric $\rho$, where

$$(3.4) \quad \rho((x_1, t_1), (x_2, t_2)) = \left| \frac{\tanh(x_1)}{1 + |t_1|} - \frac{\tanh(x_2)}{1 + |t_2|} \right| \vee |\tanh(t_1) - \tanh(t_2)|.$$

$\overline{\mathbb{R}}^2$ may be thought of as the image of $[-\infty, \infty] \times [-\infty, \infty]$ under the mapping

$$(3.5) \qquad (x, t) \rightsquigarrow (\Phi(x, t), \Psi(t)) \equiv \left( \frac{\tanh(x)}{1 + |t|}, \tanh(t) \right).$$

For $t_0 \in [-\infty, \infty]$, let $C[t_0]$ denote the set of functions $f$ from $[t_0, \infty]$ to $[-\infty, \infty]$ such that $\Phi(f(t), t)$ is continuous. Then define

$$(3.6) \qquad \Pi = \bigcup_{t_0 \in [-\infty, \infty]} C[t_0] \times \{t_0\},$$

where $(f, t_0) \in \Pi$ represents a path in $\overline{\mathbb{R}}^2$ starting at $(f(t_0), t_0)$. For $(f, t_0)$ in $\Pi$, we denote by $\hat{f}$ the function that extends $f$ to all $[-\infty, \infty]$ by setting it equal to $f(t_0)$ for $t < t_0$. Then we take

$$(3.7) \quad d((f_1, t_1), (f_2, t_2)) = \left( \sup_t |\Phi(\hat{f}_1(t), t) - \Phi(\hat{f}_2(t), t)| \right) \vee |\Psi(t_1) - \Psi(t_2)|.$$

$(\Pi, d)$ is a complete separable metric space.

Let now $\mathcal{H}$ denote the set of compact subsets of $(\Pi, d)$, with $d_\mathcal{H}$ the induced Hausdorff metric, that is,

$$(3.8) \qquad d_\mathcal{H}(K_1, K_2) = \sup_{g_1 \in K_1} \inf_{g_2 \in K_2} d(g_1, g_2) \vee \sup_{g_2 \in K_2} \inf_{g_1 \in K_1} d(g_1, g_2).$$

$(\mathcal{H}, d_\mathcal{H})$ is also a complete separable metric space.

DEFINITION 3.1. $\overline{\mathcal{W}}(\mathcal{D})$ is the closure in $(\Pi, d)$ of $\mathcal{W}(\mathcal{D})$.

Propositions 3.1 and 3.3 are essentially contained in Theorem 2.1 of [26].

PROPOSITION 3.1. $\overline{\mathcal{W}}(\mathcal{D})$ *satisfies properties* (o) *and* (i) *of Theorem* 2.1; *that is, its finite-dimensional distributions (whether from points in $\mathcal{D}$ or not) are those of coalescing Brownian motions.*

The next result is contained in Proposition B.5.

PROPOSITION 3.2. $\overline{\mathcal{W}}(\mathcal{D})$ *is almost surely a compact subset of* $(\Pi, d)$.

REMARK 3.1. Almost surely, $\overline{\mathcal{W}}(\mathcal{D}) = \lim_{k \to \infty} \mathcal{W}_k(\mathcal{D})$, where the limit is taken in $\mathcal{H}$.

REMARK 3.2. It can be shown by the methods discussed in Remark B.1 that, almost surely, all paths in $\overline{\mathcal{W}}(\mathcal{D})$ are Hölder continuous with exponent $\alpha$, for any $\alpha < \frac{1}{2}$.

PROPOSITION 3.3. *The distribution of $\overline{\mathcal{W}}(\mathcal{D})$ does not depend on $\mathcal{D}$ (including its order). Furthermore, $\overline{\mathcal{W}}(\mathcal{D})$ satisfies property* (ii) *of Theorem* 2.1.

The next theorem provides an alternative characterization to Theorem 2.1. Other characterizations that will be used for our convergence results, are presented in Section 4.

DEFINITION 3.2 (Stochastic ordering). $\mu_1 << \mu_2$ if, for $g$ any bounded measurable function on $(\mathcal{H}, \mathcal{F}_\mathcal{H})$ that is *increasing* [i.e., $g(K) \leq g(K')$ when $K \subseteq K'$], $\int g \, d\mu_1 \leq \int g \, d\mu_2$.

THEOREM 3.1. *There is an $(\mathcal{H}, \mathcal{F}_\mathcal{H})$-valued random variable $\overline{\mathcal{W}}$ whose distribution is uniquely determined by properties* (o), (i) *of Theorem* 2.1 *and*

(ii′) *if $\mathcal{W}^*$ is any other $(\mathcal{H}, \mathcal{F}_\mathcal{H})$-valued random variable satisfying* (o) *and* (i), *then $\mu_{\overline{\mathcal{W}}} \ll \mu_{\mathcal{W}^*}$.*

PROOF. The proof of this theorem follows easily from Theorem 2.1. □

**4. Characterization via counting.** In this section, we give other characterizations of the BW that will be used for our convergence theorem. They will be given in terms of the counting random variables $\eta$ and $\hat{\eta}$ defined in Definition 2.1. We begin with some properties of the BW as constructed in Section 1.

PROPOSITION 4.1. *For a BW skeleton* $\mathcal{W}(\mathcal{D})$, *the corresponding counting random variable* $\hat{\eta}_\mathcal{D} = \hat{\eta}_\mathcal{D}(t_0, t; a, b)$ *satisfies*

(4.1) $$\mathbb{P}(\hat{\eta}_\mathcal{D} \geq k) \leq \mathbb{P}(\hat{\eta}_\mathcal{D} \geq k-1)\mathbb{P}(\hat{\eta}_\mathcal{D} \geq 1)$$

(4.2) $$\leq (\mathbb{P}(\hat{\eta}_\mathcal{D} \geq 1))^k = (\Theta(b-a,t))^k,$$

*where* $\Theta(b-a,t)$ *is the probability that two independent Brownian motions starting at a distance* $b-a$ *apart at time zero will not have met by time* $t$ *(which itself can be expressed in terms of a single Brownian motion). Thus, $\hat{\eta}_\mathcal{D}$ is almost surely finite and* $\mathbb{E}(\hat{\eta}_\mathcal{D}) < \infty$.

PROOF. The proof of this proposition for $k=2$ is contained in [26]. The proof for $k>2$ can be readily obtained by using the FKG inequalities—see [12]. (The following remark notes that stronger bounds may be obtainable by the methods of [26].) □

REMARK 4.1. By analogy with the number of crossings in the scaling limit of percolation and other statistical mechanics models [2], one may expect the actual decay to be Gaussian in $k$ rather than exponential, as in (4.2). Indeed, as noted by an Associate Editor, *probably* an upper bound of the form $C_k[(b-a)/\sqrt{t}]^{k(k+1)/2}$ can be obtained by applying the method of proof of Lemma 9.4 of [26] and a result from [17].

The next proposition is a consequence of the one just before. It can also be found in Proposition 2.2 of [26].

PROPOSITION 4.2. *Almost surely, for every* $\varepsilon > 0$ *and every* $\theta = (f, t_0)$ *in* $\overline{W}(\mathcal{D})$, *there exists a path* $\theta_\varepsilon = (g, t_0')$ *in the skeleton* $\mathcal{W}(\mathcal{D})$ *such that* $g(s) = f(s)$ *for all* $s \geq t_0 + \varepsilon$.

The proof of the next proposition follows essentially from Propositions 4.1 and 4.2.

PROPOSITION 4.3. *Let* $\hat{\eta} = \hat{\eta}(t_0, t; a, b)$ *be the counting random variable for* $\overline{\mathcal{W}}(\mathcal{D})$. *Then* $\mathbb{P}(\hat{\eta} \geq k) \leq (\Theta(b-a,t))^k$, *and thus $\hat{\eta}$ is almost surely finite with finite expectation. Furthermore, $\hat{\eta} = \hat{\eta}_\mathcal{D}$ almost surely and thus*

(4.3) $$\mathbb{P}(\hat{\eta} \geq k) \leq \mathbb{P}(\hat{\eta} \geq k-1)\mathbb{P}(\hat{\eta} \geq 1)$$

(4.4) $$\leq (\mathbb{P}(\hat{\eta} \geq 1))^k = (\Theta(b-a,t))^k.$$

THEOREM 4.1. *Let $\mathcal{W}'$ be an $(\mathcal{H}, \mathcal{F}_\mathcal{H})$-valued random variable; its distribution equals that of the (standard) BW $\overline{\mathcal{W}}$ (as characterized by Theorems 2.1 and 3.1) if its finite-dimensional distributions are coalescing (standard) Brownian motions [i.e., conditions (o) and (i) of Theorem 2.1 are valid] and*

(ii″) *for all $t_0, t, a, b, \hat{\eta}_{\mathcal{W}'}$ is equidistributed with $\hat{\eta}_{\overline{\mathcal{W}}}$.*

For purposes of proving our convergence results, we will use a modified version of the above characterization theorem in which conditions (o), (i), (ii″) are all weakened.

THEOREM 4.2. *Let $\mathcal{W}'$ be an $(\mathcal{H}, \mathcal{F}_\mathcal{H})$-valued random variable and let $\mathcal{D}$ be a countable dense deterministic subset of $\mathbb{R}^2$ and for each $y \in \mathcal{D}$, let $\theta^y \in \mathcal{W}'$ be some single (random) path starting at $y$. $\mathcal{W}'$ is equidistributed with the (standard) BW $\overline{\mathcal{W}}$ if:*

(i′) *the $\theta^y$'s are distributed as coalescing (standard) Brownian motions, and*

(ii‴) *for all $t_0, t, a, b, \eta_{\mathcal{W}'} \ll \eta_{\overline{\mathcal{W}}}$, that is, $\mathbb{P}(\eta_{\mathcal{W}'} \geq k) \leq \mathbb{P}(\eta_{\overline{\mathcal{W}}} \geq k)$ for all $k$.*

PROOF. We need to show that the above conditions together imply that $\mu'$, the distribution of $\mathcal{W}'$, equals the distribution $\mu$ of the constructed BW $\overline{\mathcal{W}}$. Let $\eta'$ be the counting random variable appearing in condition (ii‴) for $\mu'$. Choose some deterministic dense countable subset $\mathcal{D}$ and consider the countable collection $\mathcal{W}^*$ of paths of $\mathcal{W}'$ starting from $\mathcal{D}$. By condition (i′), $\mathcal{W}^*$ is equidistributed with our constructed BW skeleton $\mathcal{W}$ (based on the same $\mathcal{D}$) and hence the closure $\overline{\mathcal{W}}^*$ of $\mathcal{W}^*$ in $(\Pi, d)$ is a subset of $\mathcal{W}'$ that is equidistributed with our constructed BW $\overline{\mathcal{W}}$. To complete the proof, we will use condition (ii‴) to show that $\mathcal{W}' \setminus \overline{\mathcal{W}}^*$ is almost surely empty by using the fact that the counting random variable $\eta^*$ for $\overline{\mathcal{W}}^*$ already satisfies condition (ii‴) since $\overline{\mathcal{W}}^*$ is distributed as a BW. If $\mathcal{W}' \setminus \overline{\mathcal{W}}^*$ were nonempty (with strictly positive probability), then there would have to be some rational $t_0, t, a, b$ for which $\eta' > \eta^*$. But then

$$\mathbb{P}(\eta'(t_0, t; a, b) > \eta^*(t_0, t; a, b)) > 0 \tag{4.5}$$

for some rational $t_0, t, a, b$, and this together with the fact that $\mathbb{P}(\eta' \geq \eta^*) = 1$ (which follows from $\overline{\mathcal{W}}^* \subset \mathcal{W}'$) would violate condition (ii‴) with those $t_0, t, a, b$. The proof is complete. □

REMARK 4.2. The condition $\eta_{\mathcal{W}'} \ll \eta_{\overline{\mathcal{W}}}$ can be replaced by $\mathbb{E}(\eta_{\mathcal{W}'}) \leq \mathbb{E}(\eta_{\overline{\mathcal{W}}})$. We note that $\mathbb{E}(\eta_{\overline{\mathcal{W}}}) = 1 + (b-a)/\sqrt{\pi t}$, as given in [13] by a calculation stretching back to [7]. So, in particular, in the context of Theorem 4.2, if besides (i′), $\mathbb{E}(\eta_{\mathcal{W}'}) \leq 1 + (b-a)/\sqrt{\pi t}$ for all $t_0, t, a, b$, then $\mathcal{W}'$ is equidistributed with the BW.

**5. General convergence results.** In this section, we state and prove Theorem 5.1, which is an extension of our convergence result for noncrossing paths, Theorem 2.2, to the case where paths can cross (before the scaling limit has been taken). At the end of the section, we show that the noncrossing Theorem 2.2 follows from Theorem 5.1 and other results.

Before stating our general theorem that allows crossing, we briefly discuss some systems with crossing paths, to which it should be applicable (see also Section 1.3 of [5] and [20]). Consider the stochastic difference equation (1.1) where the $\Delta_{i,j}$'s are i.i.d. integer-valued random variables, with zero mean and finite nonzero variance. Allowing $(i,j)$ to be arbitrary in $\mathbb{Z}^2$, we obtain as a natural generalization of Figure 1 a collection of random piecewise linear paths that can cross each other, but that still coalesce when they meet at a lattice point in $\mathbb{Z}^2$.

With the natural choice of diffusive space-time scaling and under conditions of irreducibility and aperiodicity (to ensure that the walks from any two starting points have a strictly positive probability of coalescing), the scaling limit of such a discrete time system should be the standard BW. To see what happens in reducible cases, consider simple random walks ($\Delta_{i,j} = \pm 1$), where the paths on the even and odd subsets of $\mathbb{Z}^2$ are independent of each other, and so the scaling limit on all of $\mathbb{Z}^2$ consists of the union of two independent BWs. For $\Delta_{i,j} = \pm 2$, the limit would be the union of four independent BWs. We remark that for continuous-time random walks (as discussed in the next section of this paper for $\Delta_{i,j} = \pm 1$), no aperiodicity condition is needed.

We proceed with some definitions needed for our general convergence theorem. For $a, b, t_0 \in \mathbb{R}$, $a < b$, and $t > 0$, we define two real-valued measurable functions $l_{t_0,t}([a,b])$ and $r_{t_0,t}([a,b])$ on $(\mathcal{H}, \mathcal{F}_\mathcal{H})$ as follows. For $K \in \mathcal{H}$, $l_{t_0,t}([a,b])$ evaluated at $K$ is defined as $\inf\{x \in [a,b] | \exists y \in \mathbb{R}$ and a path in $K$ which touches both $(x,t_0)$ and $(y,t_0+t)\}$ and $r_{t_0,t}([a,b])$ is defined similarly with the inf replaced by sup. We also define the following functions on $(\mathcal{H}, \mathcal{F}_\mathcal{H})$ whose values are subsets of $\mathbb{R}$. As before, we let $K \in \mathcal{H}$ and suppress $K$ on the left-hand side of the formula for ease of notation:

$$(5.1) \quad \mathcal{N}_{t_0,t}([a,b]) = \{y \in \mathbb{R} | \exists x \in [a,b] \text{ and a path in } K \text{ which}$$
$$\text{touches both } (x,t_0) \text{ and } (y,t_0+t)\},$$

$$(5.2) \quad \mathcal{N}_{t_0,t}^-([a,b]) = \{y \in \mathbb{R} | \text{there is a path in } K \text{ which}$$
$$\text{touches both } (l_{t_0,t}([a,b]),t_0) \text{ and } (y,t_0+t)\},$$

$$(5.3) \quad \mathcal{N}_{t_0,t}^+([a,b]) = \{y \in \mathbb{R} | \text{there is a path in } K \text{ which}$$
$$\text{touches both } (r_{t_0,t}([a,b]),t_0) \text{ and } (y,t_0+t)\}.$$

REMARK 5.1. We notice that $|\mathcal{N}_{t_0,t}([a,b])| = \eta(t_0,t;a,b)$.

Let $\{\mathcal{X}_m\}$ be a sequence of $(\mathcal{H}, \mathcal{F}_\mathcal{H})$-valued random variables with distributions $\{\mu_m\}$. We define conditions (B1'), (B2') as follows.

(B1') $\forall \beta > 0, \limsup_{m\to\infty} \sup_{t>\beta} \sup_{t_0, a\in\mathbb{R}} \mu_m(|\mathcal{N}_{t_0,t}([a-\varepsilon, a+\varepsilon])| > 1) \to 0$ as $\varepsilon \to 0^+$.

(B2') $\forall \beta > 0, \frac{1}{\varepsilon} \limsup_{m\to\infty} \sup_{t>\beta} \sup_{t_0,a\in\mathbb{R}} \mu_m(\mathcal{N}_{t_0,t}([a-\varepsilon, a+\varepsilon]) \neq \mathcal{N}^+_{t_0,t}([a-\varepsilon, a+\varepsilon]) \cup \mathcal{N}^-_{t_0,t}([a-\varepsilon, a+\varepsilon])) \to 0$ as $\varepsilon \to 0^+$.

REMARK 5.2. Note that if we consider a process with noncrossing paths, then conditions (B1') and (B2') follow from conditions (B1) and (B2), respectively, because of the following monotonicity property. For all $a < b, t_0$ and $0 < s < t$,

$$\mathbb{P}(|\eta(t_0, t; a, b)| \geq k) \leq \mathbb{P}(|\eta(t_0, s; a, b)| \geq k)$$

for all $k \in \mathbb{N}$.

THEOREM 5.1. *Suppose that $\{\mu_m\}$ is tight. If Conditions* (I1), (B1') *and* (B2') *hold, then $\{\mathcal{X}_m\}$ converges in distribution to the BW $\overline{\mathcal{W}}$.*

Theorem 5.1 is proved through a series of lemmas.

LEMMA 5.1. *Let $\mu$ be a subsequential limit of $\{\mu_m\}$ and suppose that $\mu$ satisfies condition* (i') *of Theorem* 4.2 *and:*

(B1'') $\forall \beta > 0, \sup_{t>\beta} \sup_{t_0,a} \mu(|\mathcal{N}_{t_0,t}([a-\varepsilon, a+\varepsilon])| > 1) \to 0$ *as* $\varepsilon \to 0^+$,

(B2'') $\forall \beta > 0, \frac{1}{\varepsilon} \sup_{t>\beta} \sup_{t_0,a} \mu(\mathcal{N}_{t_0,t}([a-\varepsilon, a+\varepsilon]) \neq \mathcal{N}^+_{t_0,t}([a-\varepsilon, a+\varepsilon]) \cup \mathcal{N}^-_{t_0,t}([a-\varepsilon, a+\varepsilon])) \to 0$ *as* $\varepsilon \to 0^+$.

*Then $\mu$ is the distribution of the BW.*

PROOF. It follows from conditions (i') and (B1'') that the limiting random variable $\mathcal{X}$ satisfies condition (i) of the characterization Theorem 2.1. That is, $(\mu)$ almost surely there is exactly one path starting from each point of $\mathcal{D}$ and these paths are distributed as coalescing Brownian motions. Let us define an $(\mathcal{H}, \mathcal{F}_\mathcal{H})$-valued random variable $\mathcal{X}'$ on the same probability space as the one on which $\mathcal{X}$ is defined to be the closure in $(\Pi, d)$ of the paths of $\mathcal{X}$ starting from $\mathcal{D}$. We will denote probabilities in the common probability space by $\mathbb{P}$. $\mathcal{X}'$ has the distribution of $\overline{\mathcal{W}}$. We need to show that it also satisfies condition (ii''') of Theorem 4.2. Let $a < b, t_0 \in \mathbb{R}$ and $t > 0$ be given. For the random variable $\mathcal{X}$ we will denote the counting random variable $\eta(t_0, t; a, b)$ by $\eta$ and the corresponding variable for $\mathcal{X}'$ by $\eta'$. Let $z_j = (a + j(b-a)/M, t_0)$, for $j = 0, 1, \ldots, M$, be $M+1$ equally spaced points in the interval $[a, b] \times \{t_0\}$.

Now define $\eta_M = |\{x \in \mathbb{R} | \exists$ a path in $\mathcal{X}$ which touches both a point in $\{z_0, \ldots, z_M\}$ and $(x, t+t_0)\}|$, where $|\cdot|$ stands for cardinality. Let $\eta'_M$ be the corresponding random variable for $\mathcal{X}'$. Clearly, $\eta \geq \eta_M$ and $\eta' \geq \eta'_M$. From (B1″) it follows that $\eta_M = \eta'_M$ almost surely. Now let $\varepsilon = \frac{(b-a)}{M}$. By condition (B2″), letting $M \to \infty$ ($\varepsilon \to 0$), we obtain

$$\mathbb{P}(\eta > \eta'_M) = \mathbb{P}(\eta > \eta_M) \to 0 \quad \text{as } M \to \infty.$$

Thus, $\mathbb{P}(\eta > \eta') = 0$, showing that $\eta$ is stochastically dominated by $\eta'$. This completes the proof of the lemma.

For $t > 0$, $\varepsilon > 0$, $0 < \varepsilon' < \frac{\varepsilon}{8}$, $0 \leq \delta < \frac{t}{2}$, consider the following event:

$O(a, t_0, t, \varepsilon, \varepsilon', \delta)$

$= \{K \in \mathcal{H} |$ there are three paths $(x_1(t), t_1), (x_2(t), t_2), (x_3(t), t_3)$ in $K$

with $t_1, t_2, t_3 < t_0 + \delta$, $x_1(t_0 + \delta) \in (a - \varepsilon - \varepsilon', a - \varepsilon + \varepsilon')$,

$x_2(t_0 + \delta) \in (a - \varepsilon + 2\varepsilon', a + \varepsilon - 2\varepsilon')$,

$x_3(t_0 + \delta) \in (a + \varepsilon - \varepsilon', a + \varepsilon + \varepsilon')$

and $x_2(t_0 + t) \neq x_1(t_0 + t), x_2(t_0 + t) \neq x_3(t_0 + t)\}$.

LEMMA 5.2. *Condition* (B2″) *in Lemma* 5.1 *can be replaced by:*

(B2‴) $\forall \beta > 0$, $\frac{1}{\varepsilon} \limsup_{\varepsilon' \to 0} \sup_{t > \beta} \sup_{t_0, a} \limsup_{\delta \to 0} \mu(O(a, t_0, t, \varepsilon, \varepsilon', \delta)) \to 0$ *as* $\varepsilon \to 0^+$.

PROOF. We prove the lemma by showing that conditions (i′) and (B1″) together with (B2‴) imply condition (B2″). Let $\beta > 0$. Define $C_1(b, t_0, \varepsilon', \delta)$ as

$\{K \in \mathcal{H} |$ there is a path in $K$ which touches both $(b, t_0)$

and $\{b - \varepsilon'\} \times [t_0, t_0 + \delta] \cup \{b + \varepsilon'\} \times [t_0, t_0 + \delta]\}$,

and $C_2(a, t_0, \varepsilon, \varepsilon', \delta)$ as

$\{K \in \mathcal{H} |$ there is a path in $K$ which touches both $[a - \varepsilon, a + \varepsilon] \times \{t_0\}$

and $\{a - \varepsilon - \varepsilon'\} \times [t_0, t_0 + \delta] \cup \{a + \varepsilon + \varepsilon'\} \times [t_0, t_0 + \delta]\}$.

Now observe that (modulo sets of zero $\mu$ measure)

$\{\mathcal{N}_{t_0, t}([a - \varepsilon, a + \varepsilon]) \neq \mathcal{N}^+_{t_0, t}([a - \varepsilon, a + \varepsilon]) \cup \mathcal{N}^-_{t_0, t}([a - \varepsilon, a + \varepsilon])\}$

$\cap C^c_1(a + \varepsilon, t_0, \varepsilon', \delta) \cap C^c_1(a - \varepsilon, t_0, \varepsilon', \delta) \cap C^c_2(a, t_0, \varepsilon, \varepsilon', \delta)$

$\cap \{|\mathcal{N}_{t_0 + \delta, t - \delta}([a - \varepsilon - 2\varepsilon', a - \varepsilon + 2\varepsilon'])| = 1\}$

$\cap \{|\mathcal{N}_{t_0 + \delta, t - \delta}([a + \varepsilon - 2\varepsilon', a + \varepsilon + 2\varepsilon'])| = 1\}$

$\subseteq O(a, t_0, t, \varepsilon, \varepsilon', \delta)$.

Therefore, we have

$$\mu(\mathcal{N}_{t_0,t}([a-\varepsilon, a+\varepsilon]) \neq \mathcal{N}^+_{t_0,t}([a-\varepsilon, a+\varepsilon]) \cup \mathcal{N}^-_{t_0,t}([a-\varepsilon, a+\varepsilon]))$$
$$\leq \mu(O(a,t_0,t,\varepsilon,\varepsilon',\delta)) + \mu(C_2(a,t_0,\varepsilon,\varepsilon',\delta)) + \mu(C_1(a+\varepsilon,t_0,\varepsilon',\delta))$$
$$+ \mu(C_1(a-\varepsilon,t_0,\varepsilon',\delta)) + \mu(|\mathcal{N}_{t_0+\delta,t-\delta}([a-\varepsilon-2\varepsilon', a-\varepsilon+2\varepsilon'])| > 1)$$
$$+ \mu(|\mathcal{N}_{t_0+\delta,t-\delta}([a+\varepsilon-2\varepsilon', a+\varepsilon+2\varepsilon'])| > 1).$$

Letting $\delta \to 0$, we obtain

$$\mu(\mathcal{N}_{t_0,t}([a-\varepsilon, a+\varepsilon]) \neq \mathcal{N}^+_{t_0,t}([a-\varepsilon, a+\varepsilon]) \cup \mathcal{N}^-_{t_0,t}([a-\varepsilon, a+\varepsilon]))$$
$$\leq \limsup_{\delta \to 0}\{\mu(O(a,t_0,t,\varepsilon,\varepsilon',\delta)) + \mu(C_2(a,t_0,\varepsilon,\varepsilon',\delta))$$
(5.4)
$$+ \mu(C_1(a+\varepsilon,t_0,\varepsilon',\delta)) + \mu(C_1(a-\varepsilon,t_0,\varepsilon',\delta))$$
$$+ \mu(|\mathcal{N}_{t_0+\delta,t-\delta}([a-\varepsilon-2\varepsilon', a-\varepsilon+2\varepsilon'])| > 1)$$
$$+ \mu(|\mathcal{N}_{t_0+\delta,t-\delta}([a+\varepsilon-2\varepsilon', a+\varepsilon+2\varepsilon'])| > 1)\}.$$

Now,

$$\lim_{\delta \to 0}\mu(C_1(a+\varepsilon,t_0,\varepsilon',\delta)) = \lim_{\delta \to 0}\mu(C_1(a-\varepsilon,t_0,\varepsilon',\delta))$$
$$= \lim_{\delta \to 0}\mu(C_2(a,t_0,\varepsilon,\varepsilon',\delta)) = 0,$$

since elements of $\mathcal{H}$ are compact subsets $K$ of $\Pi$, and compact sets of continuous functions are equicontinuous. If the above limit did not vanish, then there would be positive $\mu$-measure for $K$ to contain paths with arbitrarily close to flat segments, thus violating equicontinuity.

Now since, $t - \delta > \frac{t}{2} > \frac{\beta}{2}$, it follows from (B1″) that

$$\sup_{t>\beta}\sup_{a,t_0}\sup_{0<\delta<t/2} \mu(|\mathcal{N}_{t_0+\delta,t-\delta}([a-\gamma, a+\gamma])| > 1)$$
$$\leq \sup_{t>\beta/2}\sup_{a,t_0}\mu(|\mathcal{N}_{t_0,t}([a-\gamma, a+\gamma])| > 1) \to 0 \quad \text{as } \gamma \to 0.$$

This implies that for all $\varepsilon > 0$,

(5.5) $\limsup_{\varepsilon' \to 0}\sup_{t>\beta}\sup_{a,t_0}\limsup_{\delta \to 0}\mu(|\mathcal{N}_{t_0+\delta,t-\delta}([a\pm\varepsilon-2\varepsilon', a\pm\varepsilon+2\varepsilon'])| > 1) = 0.$

Together with (5.4), this gives us

$$\sup_{t>\beta}\sup_{t_0,a}\mu(\mathcal{N}_{t_0,t}([a-\varepsilon, a+\varepsilon])$$
(5.6)
$$\neq \mathcal{N}^+_{t_0,t}([a-\varepsilon, a+\varepsilon]) \cup \mathcal{N}^-_{t_0,t}([a-\varepsilon, a+\varepsilon]))$$
$$\leq \limsup_{\varepsilon' \to 0}\sup_{t>\beta}\sup_{a,t_0}\limsup_{\delta \to 0}\mu(O(a,t_0,t,\varepsilon,\varepsilon',\delta)).$$

Now, using (B2‴), we obtain

$$\frac{1}{\varepsilon}\sup_{t>\beta}\sup_{t_0,a}\mu(\mathcal{N}_{t_0,t}([a-\varepsilon,a+\varepsilon])\neq \mathcal{N}^+_{t_0,t}([a-\varepsilon,a+\varepsilon])\cup\mathcal{N}^-_{t_0,t}([a-\varepsilon,a+\varepsilon]))\to 0$$

as $\varepsilon\to 0^+$, proving the lemma. $\square$

PROOF OF THEOREM 5.1. Tightness implies that every sub-sequence of $\{\mu_m\}$ has a sub-subsequence converging to some $\mu$. Let us denote the corresponding limiting random variable by $\mathcal{X}$. We prove the theorem by showing that every such $\mu = \mu_{\overline{\mathcal{W}}}$. From Lemmas 5.1 and 5.2 it follows that it is sufficient to prove condition (i′) of Theorem 4.2, condition (B1″) and condition (B2‴).

Let $\beta > 0$ and define for all $0 \leq \delta < \frac{t}{2}$, $\mathcal{N}'^{\delta}_{t_0,t}([a,b]) = \{y \in \mathbb{R} | \exists$ a path $(x(s),s_0), s_0 < t_0+\delta$ in $K$ such that $x(t_0+\delta) \in (a,b)$ and $x(t_0+t) = y\}$. We note that the set $\{|\mathcal{N}'^{\delta}_{t_0,t}([a,b])| > 1\}$ is an open subset of $\mathcal{H}$ for all $\delta \geq 0$. Then we have

$$\sup_{t>\beta}\sup_{t_0,a}\mu(|\mathcal{N}_{t_0,t}([a-\varepsilon,a+\varepsilon])| > 1)$$

$$\leq \sup_{t>\beta}\sup_{t_0,a}\limsup_{\delta\to 0}\{\mu(|\mathcal{N}'^{\delta}_{t_0,t}([a-2\varepsilon,a+2\varepsilon])| > 1) + \mu(C_2(a,t_0,\varepsilon,\varepsilon,\delta))\}$$

$$\leq \sup_{t>\beta/2}\sup_{t_0,a}\mu(|\mathcal{N}'^{0}_{t_0,t}([a-2\varepsilon,a+2\varepsilon])| > 1)$$

$$+ \sup_{t_0,a}\limsup_{\delta\to 0}\mu(C_2(a,t_0,\varepsilon,\varepsilon,\delta)).$$

Now,

$$\limsup_{\delta\to 0}\mu(C_2(a,t_0,\varepsilon,\varepsilon,\delta)) = 0,$$

since elements of $\mathcal{H}$ are compact subsets of $\Pi$. This together with the fact that $\{|\mathcal{N}'^{\delta}_{t_0,t}([a,b])| > 1\}$ is an open subset of $\mathcal{H}$ leads to

$$\sup_{t>\beta}\sup_{t_0,a}\mu(|\mathcal{N}_{t_0,t}([a-\varepsilon,a+\varepsilon])| > 1)$$

$$\leq \sup_{t>\beta/2}\sup_{t_0,a}\mu(|\mathcal{N}'^{0}_{t_0,t}([a-2\varepsilon,a+2\varepsilon])| > 1)$$

$$\leq \sup_{t>\beta/2}\sup_{t_0,a}\limsup_{m}\mu_m(|\mathcal{N}'^{0}_{t_0,t}([a-2\varepsilon,a+2\varepsilon])| > 1)$$

$$\leq \sup_{t>\beta/2}\sup_{t_0,a}\limsup_{m}\mu_m(|\mathcal{N}_{t_0,t}([a-2\varepsilon,a+2\varepsilon])| > 1)$$

$$\leq \limsup_{m}\sup_{t>\beta/2}\sup_{t_0,a}\mu_m(|\mathcal{N}_{t_0,t}([a-2\varepsilon,a+2\varepsilon])| > 1).$$

It follows from (B1$'$) that

$$\limsup_{m} \sup_{t>\beta/2} \sup_{t_0,a} \mu_m(|\mathcal{N}_{t_0,t}([a-2\varepsilon, a+2\varepsilon])| > 1) \to 0 \qquad \text{as } \varepsilon \to 0^+.$$

This proves (B1$''$), which implies that:

(o) starting from any deterministic point, there is $\mu$-almost surely only a single path in $\mathcal{X}$.

Combining this with (I1), we readily obtain that:

(i) the finite-dimensional distributions of $\mathcal{X}$ are those of coalescing Brownian motions with unit diffusion constant.

Condition (i$'$) of Theorem 4.2 follows immediately from (o) and (i). Now we proceed to verify condition (B2$'''$).

We have

$$\sup_{t>\beta} \sup_{t_0,a} \limsup_{\delta \to 0} \mu(O(a, t_0, t, \varepsilon, \varepsilon', \delta))$$
$$\leq \sup_{t>\beta/2} \sup_{a,t_0} \mu(O(a, t_0, t, \varepsilon, \varepsilon', 0))$$
$$\leq \limsup_{m} \sup_{t>\beta/2} \sup_{a,t_0} \mu_m(O(a, t_0, t, \varepsilon, \varepsilon', 0))$$
$$\leq \limsup_{m} \sup_{t>\beta/2} \sup_{a,t_0} \mu_m(\mathcal{N}_{t_0,t}([a-\varepsilon-\varepsilon', a+\varepsilon+\varepsilon'])$$
$$\neq \mathcal{N}_{t_0,t}^+([a-\varepsilon-\varepsilon', a+\varepsilon+\varepsilon'])$$
$$\cup \mathcal{N}_{t_0,t}^-([a-\varepsilon-\varepsilon', a+\varepsilon+\varepsilon'])),$$

where the second inequality follows from the fact that $O(a, t_0, t, \varepsilon, \varepsilon', 0)$ is an open subset of $\mathcal{H}$. For the third inequality to hold we need to ensure that there is no more than one path touching either $(a-\varepsilon-\varepsilon', t_0)$ or $(a+\varepsilon+\varepsilon', t_0)$; this follows from (B1$'$). Since $\varepsilon' < \frac{\varepsilon}{8}$, $\varepsilon + \varepsilon' \to 0$ as $\varepsilon \to 0$, using condition (B2$'$), we obtain

$$\frac{1}{\varepsilon} \limsup_{\varepsilon' \to 0} \sup_{t_0,a} \limsup_{\delta \to 0} \mu(O(a, t_0, t, \varepsilon, \varepsilon', \delta)) \to 0 \qquad \text{as } \varepsilon \to 0^+,$$

proving condition (B2$'''$). This completes the proof of the theorem. $\square$

We now suppose that $\mathcal{X}_1, \mathcal{X}_2, \ldots$ is a sequence of $(\mathcal{H}, \mathcal{F}_\mathcal{H})$-valued random variables so that each $\mathcal{X}_i$ consists of *noncrossing* paths. The noncrossing condition produces a considerable simplification of Theorem 5.1, namely, Theorem 2.2.

PROOF OF THEOREM 2.2. This is an immediate consequence of Remark 5.2, Theorem 5.1 and Proposition B.2. $\square$

**6. Convergence for coalescing random walks.** We now apply Theorem 2.2 to coalescing random walks. For that, we begin by precisely defining $Y$ (resp. $\widetilde{Y}$), the set of all discrete- (resp. continuous-) time coalescing random walks on $\mathbb{Z}$. For $\delta$ an arbitrary positive real number, we obtain sets of rescaled walks, $Y^{(\delta)}$ and $\widetilde{Y}^{(\delta)}$, by the usual rescaling of space by $\delta$ and time by $\delta^2$. The (main) paths of $Y$ are the discrete-time random walks $Y_{y_0,s_0}$, as described in the Introduction and shown in Figure 1, with $(y_0, s_0) = (i_0, j_0) \in \mathbb{Z} \times \mathbb{Z}$ arbitrary except that $i_0 + j_0$ must be even. Each random walk path goes from $(i,j)$ to $(i \pm 1, j+1)$ linearly. In addition to these, we add some boundary paths so that $Y$ will be a compact subset of $\Pi$. These are all the paths of the form $(f, s_0)$ with $s_0 \in \mathbb{Z} \cup \{-\infty, \infty\}$ and $f \equiv \infty$ or $f \equiv -\infty$. Note that for $s_0 = -\infty$ there are two different paths starting from the single point at $s_0 = -\infty$ in $\overline{\mathbb{R}}^2$.

The continuous-time $\widetilde{Y}$ can be defined similarly, except that here $y_0$ is any $i_0 \in \mathbb{Z}$ and $s_0$ is arbitrary in $\mathbb{R}$. Continuous-time walks are normally seen as jumping from $i$ to $i \pm 1$ at the times $T_k^{(i)} \in (-\infty, \infty)$ of a rate-1 Poisson process. If the jump is, say, to $i+1$, then our polygonal path will have a linear segment between $(i, T_k^{(i)})$ and $(i+1, T_{k'}^{(i+1)})$, where $T_{k'}^{(i+1)}$ is the first Poisson event at $i+1$ after $T_k^{(i)}$. Furthermore, if $T_k^{(i_0)} < s_0 < T_{k+1}^{(i_0)}$, then there will be a constant segment in the path before the first nonconstant linear segment. If $s_0 = T_k^{(i_0)}$, then we take two paths: one with an initial constant segment and one without.

THEOREM 6.1. *Each of the collections of rescaled coalescing random walk paths, $Y^{(\delta)}$ (in discrete time) and $\widetilde{Y}^{(\delta)}$ (in continuous time), converges in distribution to the standard BW as $\delta \to 0$.*

PROOF. By Theorem 2.2, it suffices to verify conditions (I1), (B1) and (B2).

Condition (I1) is basically a consequence of the Donsker invariance principle, as already noted in the Introduction. Conditions (B1) and (B2) follow from the coalescing walks version of the inequality of (4.3), which is

$$(6.1) \qquad \mu_\delta(\eta(t_0, t; a, a+\varepsilon) \geq k) \leq [\mu_\delta(\eta(t_0, t; a, a+\varepsilon) \geq 2)]^{k-1}.$$

Taking the sup over $(a, t_0)$ and the lim sup over $\delta$ and using standard random walk arguments produces an upper bound of the form $C_k(\varepsilon/\sqrt{t})^{k-1}$, which yields (B1) and (B2) as desired. □

## APPENDIX A:

**Some measurability issues.** Let $(\mathcal{H}, d_\mathcal{H})$ denote the Hausdorff metric space induced by $(\Pi, d)$. $\mathcal{F}_\mathcal{H}$ denotes the $\sigma$-field generated by the open sets of $\mathcal{H}$. We will consider now *cylinders* of $\mathcal{H}$. Let us fix nonempty horizontal segments $I_1, \ldots, I_n$ in $\mathbb{R}^2$ (i.e., $I_k = I'_k \times \{t_k\}$), where each $I'_k$ is an interval (which need not be finite and can be open, closed or neither) and $t_k \in \mathbb{R}$. Define

$$\begin{aligned}
C^{t_0}_{I_1,\ldots,I_n} := \{K \in \mathcal{H}\colon \text{ there exists } (f,t) \in K \text{ such that} \\
t > t_0 \text{ and } (f,t) \text{ goes through } I_1, \ldots, I_n\},
\end{aligned} \tag{A.1}$$

$$\begin{aligned}
\overline{C}^{t_0}_{I_1,\ldots,I_n} := \{K \in \mathcal{H}\colon \text{ there exists } (f,t) \in K \text{ such that} \\
t \geq t_0 \text{ and } (f,t) \text{ goes through } I_1, \ldots, I_n\},
\end{aligned} \tag{A.2}$$

$$\begin{aligned}
C_{I_1,\ldots,I_n} := \{K \in \mathcal{H}\colon \text{ there exists } (f,t) \in K \text{ such that} \\
(f,t) \text{ goes through } I_1, \ldots, I_n\}.
\end{aligned} \tag{A.3}$$

We will call sets of the form (A.1) *open cylinders* if each $I_k$ is open, and sets of the form (A.2) *closed cylinders* if each $I_k$ is closed.

REMARK A.1. It is easy to see that sets of the form (A.1)–(A.3) for arbitrary $I_1, \ldots, I_n$ can be generated by open cylinders.

Let now $\mathfrak{C}$ be the $\sigma$-field generated by the open cylinders.

PROPOSITION A.1. $\mathcal{F}_\mathcal{H} = \mathfrak{C}$.

The proposition is a consequence of the following two lemmas.

LEMMA A.1. $\mathcal{F}_\mathcal{H} \supset \mathfrak{C}$.

PROOF. It is enough to observe that the open cylinders are open sets of $\mathcal{H}$. Indeed, take an open cylinder, an element $K$ in that cylinder, and $(f,t) \in K$ such that $t > t_0$ and $a_i < f(t_i) < b_i$ for all $i = 1, \ldots, n$. All points of $B_\mathcal{H}(K, \varepsilon)$, the open ball in $\mathcal{H}$ around $K$ with radius $\varepsilon$, contain a path $(f', t')$ in a ball in $\Pi$ around $(f,t)$ of radius $\varepsilon$. Thus by choosing $\varepsilon$ small enough, $(f', t')$ will satisfy $t_0 < t' < t_i$ and $a_i < f'(t_i) < b_i$ for all $i = 1, \ldots, n$. □

LEMMA A.2. $\mathcal{F}_\mathcal{H} \subset \mathfrak{C}$.

PROOF. It is enough to generate the $\varepsilon$-balls in $\mathcal{H}$ with cylinders. We will start with $\varepsilon$-balls around points of $\mathcal{H}$ consisting of finitely many paths of $\Pi$.

We will use the concept of a *cone* in $\mathbb{R}^2$ around $(f,t)$. Let $r^- = r^-(t,\varepsilon)$ and $r^+ = r^+(t,\varepsilon)$ be the two solutions of

$$|\tanh(r) - \tanh(t)| = \varepsilon, \tag{A.4}$$

with $r^- \leq r^+$. For $s$ fixed, let $x^-(s) = x^-(s, \varepsilon)$ and $x^+(s) = x^+(s, \varepsilon)$ be the solutions for small $\varepsilon$ of

$$\text{(A.5)} \qquad \frac{\tanh(x) - \tanh(\hat{f}(s))}{|s| + 1} = \pm \varepsilon,$$

with $x^-(s) \leq x^+(s)$. The cone around $(f, t)$ is defined as

$$\text{(A.6)} \qquad \mathbb{C} := \{(x, y) \in \mathbb{R}^2 : x^-(y) \leq x \leq x^+(y), y \geq r^-\}.$$

Now let $K_0 = \{(f_1, t_1), \ldots, (f_n, t_n)\}$. Let $\mathbb{C}_1, \ldots, \mathbb{C}_n$ be the respective cones of $(f_1, t_1), \ldots, (f_n, t_n)$. For $i = 1, \ldots, n$, let $r_i^+ = r^+(t_i, \varepsilon)$, $r_i^- = r^-(t_i, \varepsilon)$.

Consider now a family of horizontal lines $\{L_1, L_2, \ldots\} = \mathbb{R} \times \mathcal{S}$, where $\mathcal{S} = \{s_1, s_2, \ldots\}$, with the $s_k$'s distinct and such that $\bigcup_{k \geq 1} L_k$ is dense in $\mathbb{R}^2$. For fixed $k$, consider the segments (of nonzero length) $I_k^i$ into which $L_k$ is divided by all the points of the form $x_i^-(s_k)$ and $x_i^+(s_k)$ for $i = 1, \ldots, n$ (number of such segments $\leq 2n + 1$). Segments with interior points in some cone are closed; otherwise, they are open.

Let $\mathcal{I} = \{I_{k_1}^{i_1}, \ldots, I_{k_m}^{i_m}\}$ be any finite sequence of the intervals defined above, with $k_1, \ldots, k_m$ distinct. For $1 \leq i \leq n$, we will say that $\mathcal{I}$ is $i$-good if $\mathbb{C}_i$ contains all the intervals in $\mathcal{I}$. If $\mathcal{I}$ is not $i$-good for any $1 \leq i \leq n$, then $\mathcal{I}$ is $bad$. [Roughly speaking, $\mathcal{I}$ is bad unless the intervals of $\mathcal{I}$ closely track (within distance $\varepsilon$) some particular path of $K_0$.] Let

$$\text{(A.7)} \quad \widehat{C}_i := \{K \in \mathcal{H} : \text{there exists } (f, t) \in K \text{ such that}$$
$$t \in [r_i^-, r_i^+] \text{ and } \{f(s)\} \times \{s\} \in \mathbb{C}_i \text{ for all } s \geq t\}$$

$$\text{(A.8)} \quad = \{K \in \mathcal{H} : \text{there exists } (f, t) \in K \text{ such that } d((f, t), (f_i, t_i)) \leq \varepsilon\}.$$

It is not hard to see that $\widehat{C}_i$ belongs to $\mathfrak{C}$ by writing $[r_i^-, r_i^+]$ as a finite union of small subintervals and then approximating $\widehat{C}_i$ by a finite union of sets of the form $\overline{C}_\mathcal{I}^s$, where $\mathcal{I} = \{I_1, \ldots, I_m\}$; $I_\ell = [x_i^-(s_\ell'), x_i^+(s_\ell')] \times \{s_\ell'\}$; $s \leq s_1' \leq s_2' \leq \cdots$; $s, s_1' \in [r_i^-, r_i^+]$; and each $s_\ell' \in \mathcal{S}$. Note that in the definition (A.2) of such a $\overline{C}_\mathcal{I}^s$, the starting time $t$ of the path $(f, t)$ must be in $[s, s_1']$. Define $\widehat{C} := \bigcap_{i=1}^n \widehat{C}_i$.

We next give an explicit, somewhat complicated, formula for the closed ball $\overline{B}_\mathcal{H}(K_0, \varepsilon)$. An explanation is presented immediately after the formula:

$$\text{(A.9)} \quad \overline{B}_\mathcal{H}(K_0, \varepsilon) = \left[\left(\bigcup_{m, \mathcal{I} : \mathcal{I} \text{ is bad}} C_\mathcal{I}\right) \cup \left(\bigcup_{i=1}^n \bigcup_{\substack{m, \mathcal{I}, k : \mathcal{I} \text{ is } i\text{-good and} \\ s_k > \max_j \{r_j^+ : \mathcal{I} \text{ is } j\text{-good}\}}} C_\mathcal{I}^{s_k}\right)\right]^c \cap \widehat{C}.$$

If $K \in \mathcal{H}$ is such that $d_\mathcal{H}(K, K_0) \leq \varepsilon$, then (i) for each $i = 1, \ldots, n$, there exists a path $(f, t) \in K$ such that $d((f, t), (f_i, t_i)) \leq \varepsilon$, which is clearly equivalent to $K$ belonging to $\widehat{C}$, and (ii) for each path $(f, t) \in K$, there exists $i \in \{1,\ldots,n\}$ such that $d((f,t),(f_i,t_i)) \leq \varepsilon$. The latter condition is equivalent to (ii′) there is no path in $K$ which is at a distance greater than $\varepsilon$ from every $(f_i,t_i)$, $i=1,\ldots,n$. But that is equivalent to $K$ not belonging to the set within square brackets. Indeed, bad $\mathcal{I}$'s, in the first term of that expression, ensure that some path in $K$ is at distance greater than $\varepsilon$ of $(f_i,t_i)$ for every $i=1,\ldots,n$ spatially; in the second term, some path in $K$ starts at a distance greater than $\varepsilon$ from the starting time of every $(f_j,t_j)$ from which its spatial distance is acceptable. Equation (A.9) is thus established.

To complete the proof we generalize from $K_0$ finite to a general $K$. To generate $B_{\mathcal{H}}(K,\varepsilon)$ for arbitrary $K \in \mathcal{H}$, we approximate $B_{\mathcal{H}}(K,\varepsilon)$ by an increasing sequence of balls around $\widetilde{K}$'s consisting of finitely many paths. For that, we note that, by compactness of $K$, for every integer $j > 1$, there exists $K_j \in \mathcal{H}$ consisting of finitely many paths such that $K_j \subset K$ (as subsets of $\Pi$) and $d_{\mathcal{H}}(K,K_j) < \varepsilon/j$ for all $j > 1$. We then have

$$\text{(A.10)} \qquad B_{\mathcal{H}}(K,(1-2/j)\varepsilon) \subset B_{\mathcal{H}}(K_j,(1-1/j)\varepsilon) \subset B_{\mathcal{H}}(K,\varepsilon).$$

The first inclusion is justified as follows. Let $K' \in B_{\mathcal{H}}(K,(1-1/j)\varepsilon)$. Then $d_{\mathcal{H}}(K,K') < (1-1/j)\varepsilon$ and, by the triangle inequality,

$$\text{(A.11)} \qquad \begin{aligned} d_{\mathcal{H}}(K_j,K') &\leq d_{\mathcal{H}}(K_j,K) + d_{\mathcal{H}}(K,K') \\ &< \varepsilon/j + (1-2/j)\varepsilon = (1-1/j)\varepsilon. \end{aligned}$$

Thus $K' \in B_{\mathcal{H}}(K_j,\varepsilon)$. The second inclusion is justified similarly. It is clear now that $\bigcup_{j>1} B_{\mathcal{H}}(K,(1-2/j)\varepsilon) = \bigcup_{j>1} B_{\mathcal{H}}(K_j,(1-1/j)\varepsilon) = B_{\mathcal{H}}(K,\varepsilon)$. □

## APPENDIX B:

**Compactness and tightness.** Let $\Lambda_{L,T} = [-L,L] \times [-T,T]$, and let $\{\mu_m\}$ be a sequence of probability measures on $(\mathcal{H},\mathcal{F}_{\mathcal{H}})$. For $x_0,t_0 \in \mathbb{R}$ and $u,t > 0$, let $R(x_0,t_0;u,t)$ denote the rectangle $[x_0 - \frac{u}{2}, x_0 + \frac{u}{2}] \times [t_0, t_0 + t]$ in $\mathbb{R}^2$. Define $A_{t,u}(x_0,t_0)$ to be the event (in $\mathcal{F}_{\mathcal{H}}$) that $K$ (in $\mathcal{H}$) contains a path touching both $R(x_0,t_0;\frac{u}{2},t)$ and (at a later time) the left or right boundary of the bigger rectangle $R(x_0,t_0;u,2t)$. See Figure 2.

Our tightness condition is

(T1) $\tilde{g}(t,u;L,T) \equiv t^{-1} \limsup_m \sup_{(x_0,t_0) \in \Lambda_{L,T}} \mu_m(A_{t,u}(x_0,t_0)) \to 0$ as $t \to 0^+$.

PROPOSITION B.1. *Condition* (T1) *implies tightness of* $\{\mu_m\}$.

PROOF. Let

$$g_m(t,u;L,T) = \sup_{(x_0,t_0) \in \Lambda_{L,T}} \mu_m(A_{t,u}(x_0,t_0)).$$

Now define $B_{t,u}(x_0, t_0)$ as the event (in $\mathcal{F}_\mathcal{H}$) that $K$ (in $\mathcal{H}$) contains a path which touches a point $(x', t') = (f(t'), t') \in R(x_0, t_0; \frac{u}{2}, t)$ and for some $t'' \in [t', t' + t], |f(t'') - f(t')| \geq u$. We observe that $B_{t,u}(x_0, t_0) \subseteq A_{t,u}(x_0, t_0)$.

We now cover $\Lambda_{L,T}$ with $\frac{u}{2} \times t$ rectangular boxes. Let $L_D = L_D(u) = \{-L + k\frac{u}{2} : k \in \mathbb{Z}, 0 \leq k \leq \lceil \frac{2L}{u/2} \rceil\}$ and $T_D = T_D(t) = \{-T + mt : m \in \mathbb{Z}, 0 \leq m \leq \lceil \frac{2T}{t} \rceil\}$. Then,

$$\mu_m\left(\bigcup_{(x_0,t_0)\in\Lambda_{L,T}} B_{t,u}(x_0,t_0)\right)$$

$$\leq \mu_m\left(\bigcup_{(x_0,t_0)\in L_D\times T_D} B_{t,u}(x_0,t_0)\right)$$

$$\leq \mu_m\left(\bigcup_{(x_0,t_0)\in L_D\times T_D} A_{t,u}(x_0,t_0)\right)$$

(B.1)
$$\leq \left\lceil\frac{2L+1}{u/2}\right\rceil\left\lceil\frac{2T+1}{t}\right\rceil g_m(t,u;L,T)$$

$$\leq C'\frac{LT}{tu}g_m(t,u;L,T)$$

$$\leq C'\frac{LT}{u}(\tilde{g}(t,u;L,T) + \delta)$$

for any $\delta > 0$, where in (B.1) $m$ is larger than some $M(t, u; L, T; \delta)$. The first inequality follows from the observation that if $K$ (in $\mathcal{H}$) is an outcome in $B_{t,u}(x,t)$ for some $(x,t) \in \Lambda_{L,T}$, then $K$ is an outcome in $B_{t,u}(x',t')$ for some $(x', t') \in L_D \times T_D$.

The strategy of the remainder of the proof is to use (B.1) to control the oscillations of paths within a single large rectangle $\Lambda_{L,T}$ and then, by the compactification of $\mathbb{R}^2$ [see (3.4) and (3.5)], we will control the oscillations globally by appropriately choosing sequences of $L, T$ values tend to and $t, u, \delta$ values tend to 0.

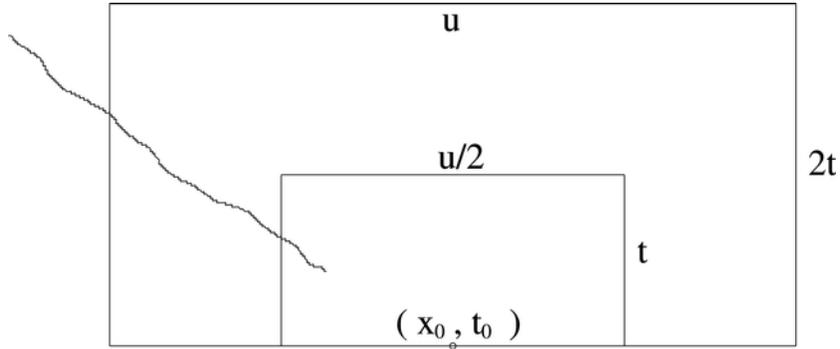

Fig. 2. *Schematic diagram of a path causing the unlikely event $A_{t,u}(x_0, t_0)$ to occur.*

Now let $\{u_n\}$ be a sequence of positive real numbers with $\lim_{n \to \infty} u_n = 0$. Since $\Phi(x,t) = (1+|t|)^{-1} \tanh(x)$, it easily follows that we can choose $L_n \to \infty$ and $T_n \to \infty$ such that $|\Phi(x,t)| \leq u_n$ if $|t| \geq T_n$ and $|\Phi(x,t) - \Phi(\pm L_n, t)| \leq u_n$ if $|t| < T_n$ and $\pm x \geq L_n$. Now choose sequences of positive real numbers $\{t'_n\}, \{\delta_n\} \to 0$ such that $C' \frac{L_n T_n}{u_n}(\tilde{g}(t'_n, u_n, L_n, T_n) + \delta_n) \leq 2^{-n}$. For all $n \in \mathbb{N}$, let

$$(B.2) \quad C_n(t) = C(t, u_n; L_n, T_n) \equiv \bigcup_{(x_0, t_0) \in \Lambda_{L_n, T_n}} B_{t, u_n}(x_0, t_0).$$

From (B.1) we have $\mu_m(C_n(t'_n)) \leq 2^{-n}$ if $m \geq M_n := M(t'_n, u_n; L_n, T_n; \delta_n)$. Let $\{t''_k\}$ be a sequence of real numbers converging to 0. If $\bigcap_{k=1}^{\infty} C_n(t''_k) \neq \varnothing$, then there exists a compact subset $K''$ of $\Pi$ which belongs to $C_n(t''_k)$ for all $k$. Since a compact set of continuous functions is equicontinuous, this is impossible and we conclude that $\bigcap_{k=1}^{\infty} C_n(t''_k) = \varnothing$. Therefore, $\mu_m(C(t, u; L, T)) \to 0$ as $t \to 0$ for any fixed $m, u, L, T$. Thus there exists $t''_n > 0$ such that $\mu_m(C_n(t''_n)) \leq 2^{-n}$ for all $m \leq M_n$. If we now let $C_n = C_n(t'_n \wedge t''_n)$, then, by the monotonicity of $C_n(t)$ in $t$, we have for all $m$,

$$(B.3) \quad \mu_m(C_n) \leq \mu_m(C_n(t'_n)) \wedge \mu_m(C_n(t''_n)) \leq 2^{-n}.$$

Now let $K \in C_n^c$ be a compact set of paths. A bound on the oscillations of paths in $K$ can be obtained as follows. Let $\psi_n = \Psi(T_n + t_n) - \Psi(T_n)$. [Recall that $\Psi(t) = \tanh(t)$.] Suppose $(f, t_0) \in K$. If $t_0 \leq s_1 \leq s_2$ are times such that $|\Psi(s_2) - \Psi(s_1)| \leq \psi_n$, then $|\Phi(f(s_1), s_1) - \Phi(f(s_2), s_2)| \leq 3u_n$. [E.g., note that $|\Psi(s_2) - \Psi(s_1)| \leq \psi_n$ for $|s_1|, |s_2| \leq T_n$ implies $|s_2 - s_1| \leq t_n$.] Let $G_n = \bigcap_{i=n+1}^{\infty} C_i^c$. Then for any $m$,

$$(B.4) \quad \mu_m(G_n) = 1 - \mu_m\left(\bigcup_{i=n+1}^{\infty} C_i\right) \geq 1 - \sum_{i=n+1}^{\infty} 2^{-i} = 1 - 2^{-n}.$$

Finally, let $D_n = \bigcup_{K \in G_n} K$. Then $D_n$ is a family of equicontinuous functions. By the Arzela–Ascoli theorem, $D_n$ is a compact subset of $\Pi$. Since $G_n$ is a collection of closed (and hence compact) subsets of $D_n$, $G_n$ is a compact subset of $\mathcal{H}$. Let $\varepsilon > 0$. Choose $n(\varepsilon) \in \mathbb{N}$ such that $2^{-n(\varepsilon)} \geq \varepsilon$. Then we have

$$(B.5) \quad \sup_m \mu_m(G_{n(\varepsilon)}) \geq 1 - \varepsilon,$$

where $G_{n(\varepsilon)}$ is a compact subset of $\mathcal{H}$. This proves that the family of measures $\{\mu_m\}$ is tight. $\square$

REMARK B.1. An argument similar to that for Proposition B.1 can be made to show that, if instead of (T1), one has the condition

(T1′) $\sum_{t\,:\,t=2^{-k},k\in\mathbb{N}} t^{-(1+\alpha)} \sup_m \sup_{x_0,t_0} \mu_m(A_{t,t^\alpha}(x_0,t_0)) < \infty$

for some $\alpha > 0$, then each $\mu_m$ as well as any subsequential limit $\mu$ of $(\mu_m)$ is supported on paths which are Hölder continuous with index $\alpha$.

PROPOSITION B.2. *Suppose $\{\mathcal{X}_m\}$ is a sequence of $(\mathcal{H}, d_\mathcal{H})$-valued random variables whose paths are noncrossing. Suppose in addition,*

(I1′) *For each $y \in \mathcal{D}$, there exist (measurable) path-valued random variables $\theta_m^y \in \mathcal{X}_m$ such that $\theta_m^y$ converges in distribution to a Brownian motion $Z_y$ starting at $y$.*

*Then the distributions $\{\mu_m\}$ of $\{\mathcal{X}_m\}$ are tight.*

PROOF. From the proof of Proposition B.1, it is sufficient to show that for each $u > 0$,

$$\limsup_m \mu_m\left(\bigcup_{(x_0,t_0) \in L_D(u) \times T_D(t)} B_{t,u}(x_0,t_0)\right) \to 0 \quad \text{as } t \to 0.$$

For $u > 0$, $t > 0$, $(x_0, t_0) \in \mathbb{R}^2$, choose two points $y_1, y_2 \in \mathcal{D}$ from the two rectangles $R(x_0 \mp \frac{3}{8}u, t_0 - \frac{t}{2}; \frac{u}{8}, \frac{t}{4})$, respectively. Let

$$B_1^m(x_0, t_0, t, u) = \left\{K \in \mathcal{H} \,\Big|\, \max_{s \le t_0 + 2t} |\theta_m^{y_1}(s) - y_1| < \frac{u}{16}\right\},$$

$$B_2^m(x_0, t_0, t, u) = \left\{K \in \mathcal{H} \,\Big|\, \max_{s \le t_0 + 2t} |\theta_m^{y_2}(s) - y_2| < \frac{u}{16}\right\},$$

and $D_{t,u}^m(x_0, t_0) = B_1^m \cap B_2^m$. Now observe that $D_{t,u}^m(x_0, t_0) \subseteq B_{t,u}^c(x_0, t_0)$ for large enough $m$. Therefore we have

$$\text{(B.6)} \qquad \limsup_m \mu_m\left(\bigcup_{(x_0,t_0) \in L_D \times T_D} B_{t,u}(x_0,t_0)\right)$$

$$\text{(B.7)} \qquad \le \sum_{(x_0,t_0) \in L_D \times T_D} \left[1 - \liminf_m \mu_m(D_{t,u}^m(x_0,t_0))\right].$$

Since $\theta_m^y$ converges in distribution to a Brownian motion $Z_y$ starting at $y$, we have

$$\text{(B.8)} \qquad \liminf_m (\mu_m(B_1^m)) = \mathbb{P}\left(\max_{s \le t_0 + 2t} |Z_{y_1}(s) - y_1| < \frac{u}{16}\right)$$

$$\text{(B.9)} \qquad \ge 1 - \frac{Ct^2}{u^4}$$

and

(B.10) $$\liminf_m(\mu_m(B_2^m)) = \mathbb{P}\bigg(\max_{s \leq t_0 + 2t} |Z_{y_2}(s) - y_2| < \frac{u}{16}\bigg)$$

(B.11) $$\geq 1 - \frac{Ct^2}{u^4}.$$

Therefore we have

$$\liminf_m \mu_m(D_{t,u}^m(x_0, t_0)) \geq 1 - 2Ct^2/u^4,$$

which gives us

$$\limsup_m \mu_m\bigg(\bigcup_{(x_0,t_0) \in L_D(u) \times T_D(t)} B_{t,u}(x_0, t_0)\bigg) \leq 2C \sum_{(x_0,t_0) \in L_D(u) \times T_D(t)} t^2/u^4.$$

Since $|L_D(u) \times T_D(t)| \sim \frac{1}{ut}$, we have shown that

$$\limsup_m \mu_m\bigg(\bigcup_{(x_0,t_0) \in L_D(u) \times T_D(t)} B_{t,u}(x_0, t_0)\bigg) \to 0 \qquad \text{as } t \to 0,$$

and the proof is complete. $\square$

REMARK B.2. The proof of Proposition B.2 shows that the limiting processes $Z_y$ starting at $y = (\bar{x}, \bar{t})$ need not be Brownian motions. It is sufficient that they be continuous processes such that for each fixed $u > 0$,

(B.12) $$\frac{1}{t} \sup_y \mathbb{P}\bigg(\sup_{\bar{t} \leq s \leq \bar{t}+t} |Z_y(s) - Z_y(\bar{t})| \geq u\bigg) \to 0 \qquad \text{as } t \to 0^+.$$

PROPOSITION B.3. *Let $\mathcal{D}$ be a countable dense subset of $\mathbb{R}^2$ and let $\mu_k$ be the distribution of the $(\mathcal{H}, \mathcal{F}_\mathcal{H})$-valued random variable $\mathcal{W}_k = \mathcal{W}_k(\mathcal{D}) = \{\widetilde{\mathcal{W}}_1, \ldots, \widetilde{\mathcal{W}}_k\}$ [as defined in (3.2)]. Then the family of measures $\{\mu_k\}$ is tight.*

PROOF. This is an immediate consequence of Proposition B.2. $\square$

PROPOSITION B.4. *If $\mathcal{W}_n$ is an a.s. increasing sequence of $(\mathcal{H}, d_\mathcal{H})$-valued random variables and the family of distributions $\{\mu_n\}$ of $\mathcal{W}_n$ is tight, then $\overline{\bigcup_n \mathcal{W}_n}$ is almost surely compact [in $(\Pi, d)$].*

PROOF. Let $\widetilde{\mathcal{W}}_k$ be an increasing sequence of points (subsets of $\Pi$) in $(\mathcal{H}, d_\mathcal{H})$, which converge in $d_\mathcal{H}$ metric to some point $\widetilde{\mathcal{W}}$ in $(\mathcal{H}, d_\mathcal{H})$. If for some $k$, $\widetilde{\mathcal{W}}_k$ is not a subset of $\widetilde{\mathcal{W}}$, then there exists an $\varepsilon > 0$ such that $d_\mathcal{H}(\widetilde{\mathcal{W}}, \widetilde{\mathcal{W}}_n) > \varepsilon$ for all $n \geq k$, contradicting the claim that $\widetilde{\mathcal{W}}_k$ converges to $\widetilde{\mathcal{W}}$. Therefore, $\widetilde{\mathcal{W}}_k \subseteq \widetilde{\mathcal{W}}$ for all $k$. This implies $\overline{\bigcup_k \widetilde{\mathcal{W}}_k} \subseteq \widetilde{\mathcal{W}}$ and therefore is a compact subset of $\Pi$ since it is a closed subset of the compact set $\widetilde{\mathcal{W}}$. Since $\{\mu_n\}$ is tight, given an $\varepsilon > 0$, there exists a compact subset K of $\mathcal{H}$ such that $\mathbb{P}(\mathcal{W}_n \in K) \geq 1 - \varepsilon$ for all $n$, so by monotonicity, $P(\mathcal{W}_n \in K$ for all $n) \geq 1 - \varepsilon$. But if $\mathcal{W}_n \in K$ for all $n$, then since $K$ is compact, there exists a subsequence $\mathcal{W}_{n_j}$ which converges to a point in $K$ and thus in $\mathcal{H}$. This implies by the first part of this proof that $\overline{\bigcup_{n_j} \mathcal{W}_{n_j}}$ $(= \overline{\bigcup_n \mathcal{W}_n}$ because $\mathcal{W}_n$ is increasing in $n$) is a compact subset of $\Pi$. Thus we have shown that $\mathbb{P}(\overline{\bigcup_n \mathcal{W}_n}$ is a compact subset of $\Pi) \geq 1 - \varepsilon$. Since the claim is true for all $\varepsilon > 0$, we have proved the proposition. □

PROPOSITION B.5. *Let $\widehat{\mathcal{D}} = \{(\hat{x}_i, \hat{t}_i) : i = 1, 2, \ldots\}$ be a (deterministic) dense countable subset of $\mathbb{R}^2$ and let $\{\widehat{\mathcal{W}}_i : i = 1, 2, \ldots\}$ be $(\Pi, d)$-valued random variables starting from $(\hat{x}_i, \hat{t}_i)$. Suppose that the joint distribution of each finite subset of the $\widehat{\mathcal{W}}_i$'s is that of coalescing Brownian motions. Then $\overline{\bigcup_{n=1}^\infty \{\widehat{\mathcal{W}}_1, \widehat{\mathcal{W}}_2, \ldots, \widehat{\mathcal{W}}_n\}}$ is almost surely compact. In particular, for $\mathcal{W}_n$ defined in Proposition (B.3), $\overline{\mathcal{W}} = \overline{\bigcup_n \mathcal{W}_n}$ is almost surely compact.*

PROOF. The proof follows immediately from Propositions B.3 and B.4. □

PROPOSITION B.6. *Let $\widehat{\mathcal{D}}$ and $\{\widehat{\mathcal{W}}_i : i = 1, 2, \ldots\}$ be as in Proposition B.5 and let $\{\widehat{\mathcal{W}}'_i : i = 1, 2, \ldots\}$ (on some other probability space) be equidistributed with $\{\widehat{\mathcal{W}}_i : i = 1, 2, \ldots\}$. Then $\widehat{\mathcal{W}} \equiv \overline{\{\widehat{\mathcal{W}}_i : i = 1, 2, \ldots\}}$ and $\widehat{\mathcal{W}}' \equiv \overline{\{\widehat{\mathcal{W}}'_i : i = 1, 2, \ldots\}}$ are equidistributed $(\mathcal{H}, \mathcal{F}_\mathcal{H})$-valued random variables.*

PROOF. It is an easy consequence of Proposition B.5 that $\{\widehat{\mathcal{W}}_i : i = 1, \ldots, n\}$ (resp. $\{\widehat{\mathcal{W}}'_i : i = 1, \ldots, n\}$) converges a.s. as $n \to \infty$ in $(\mathcal{H}, d_\mathcal{H})$ to $\widehat{\mathcal{W}}$ (resp. $\widehat{\mathcal{W}}'$). But then the identical distributions of $\{\widehat{\mathcal{W}}_i : i = 1, \ldots, n\}$ and $\{\widehat{\mathcal{W}}'_i : i = 1, \ldots, n\}$ converge, respectively, to the distributions of $\widehat{\mathcal{W}}$ and $\widehat{\mathcal{W}}'$, which thus must be identical. □

**Acknowledgments.** The authors thank Raghu Varadhan for many useful discussions, Balint Tóth for comments related to [26] and Rahul Roy, Anish Sarkar, Ed Waymire and Larry Winter for discussions and references on drainage networks. They thank Rongfeng Sun for pointing out an error in an ealier version of Appendix B. They also thank the Abdus Salam International Centre for Theoretical Physics for its hospitality when an early draft of this work was being prepared. L. R. G. Fontes thanks the Courant Institute of NYU and the Math Department at the University of Rome "La Sapienza" for hospitality and support during visits where parts of this work were done. M. Isopi thanks the Courant Institute of NYU and the IME at the University of São Paulo for hospitality and support during visits where parts of this work were done. K. Ravishankar thanks the Courant Institute of NYU, where parts of this work were done, for hospitality and support during a sabbatical leave visit.

## REFERENCES

[1] Aizenman, M. (1998). Scaling limit for the incipient spanning clusters. In *Mathematics of Multiscale Materials* (K. M. Golden, G. R. Grimmett, R. D. James, G. W. Milton and P. N. Sen, eds.) 1–24. Springer, New York. MR1635999

[2] Aizenman, M. and Burchard, A. (1999). Hölder regularity and dimension bounds for random curves. *Duke Math. J.* **99** 419–453. MR1712629

[3] Aizenman, M., Burchard, A., Newman, C. M. and Wilson, D. B. (1999). Scaling limits for minimal and random spanning trees in two dimensions. *Random Structures Algorithms* **15** 319–367. MR1716768

[4] Arratia, R. (1979). Coalescing Brownian motions on the line. Ph.D. dissertation, Univ. Wisconsin, Madison.

[5] Arratia, R. (1981). Coalescing Brownian motions and the voter model on $\mathbb{Z}$. Unpublished partial manuscript. Available from rarratia@math.usc.edu.

[6] Arratia, R. (1981). Limiting point processes for rescalings of coalescing and annihilating random walks on $\mathbb{Z}^d$. *Ann. Probab.* **9** 909–936. MR632966

[7] Bramson, M. and Griffeath, D. (1980). Clustering and dispersion rates for some interacting particle systems on $\mathbb{Z}^1$. *Ann. Probab.* **8** 183–213. MR566588

[8] Burdzy, K. and Le Gall, J.-F. (2001). Super-Brownian motion with reflecting historical paths. *Probab. Theory Related Fields* **121** 447–491. MR1872425

[9] Burdzy, K. and Mytnik, L. (2002). Super-Brownian motion with reflecting historical paths. II Convergence of approximations. Preprint. Available at www.math.washington.edu/~burdzy/preprints.shtml. MR1932696

[10] Donsker, M. D. (1951). An invariance principle for certain probability limit theorems. *Mem. Amer. Math. Soc.* **6** 1–12. MR40613

[11] Ferrari, P. A., Fontes, L. R. G. and Wu, X. Y. (2003). Two dimensional Poisson trees converge to the Brownian web. Available at arxiv.org/abs/math.PR/0304247.

[12] Fontes, L. R. G., Isopi, M., Newman, C. M. and Ravishankar, K. (2003). The Brownian web: Characterization and convergence. Available at arxiv.org/abs/math.PR/0304119.

[13] Fontes, L. R. G., Isopi, M., Newman, C. M. and Ravishankar, K. (2002). The Brownian web. *Proc. Nat. Acad. Sci. U.S.A.* **99** 15888–15893. MR1944976

[14] Fontes, L. R. G., Isopi, M., Newman, C. M. and Stein, D. L. (2001). 1D aging. Unpulished manuscript.

[15] Fontes, L. R. G., Isopi, M., Newman, C. M. and Stein, D. L. (2001). Aging in 1D discrete spin models and equivalent systems. *Phys. Rev. Lett.* **87** 110201-1–110201-4.

[16] Gangopadhyay, S., Roy, R. and Sarkar, A. (2004). Random oriented trees: A model of drainage networks. *Ann. Appl. Probab.* **14** 1242–1266. MR2071422

[17] Grabiner, D. (1999). Brownian motion in a Weyl chamber, non-colliding particles, and random matrices. *Ann. Inst. H. Poincaré Probab. Statist.* **35** 177–204. MR1678525

[18] HARRIS, T. E. (1978). Additive set-valued Markov processes and graphical methods. *Ann. Probab.* **6** 355–378. MR488377

[19] HARRIS, T. E. (1984). Coalescing and noncoalescing flows in $\mathbb{R}^1$. *Stochastic Process. Appl.* **17** 187–210. MR751202

[20] NEWMAN, C. M., RAVISHANKAR, K. and SUN, R. (2004). Convergence of coalescing nonsimple random walks to the Brownian web. Unpublished manuscript.

[21] NGUYEN, B. (1990). Percolation of coalescing random walks. *J. Appl. Probab.* **27** 269–277.

[22] PITERBARG, V. V. (1998). Expansions and contractions of isotropic stochastic flows of homeomorphisms. *Ann. Probab.* **26** 479–499. MR1626162

[23] RODRIGUEZ-ITURBE, I. and RINALDO, A. (1997). *Fractal River Basins*: *Chance and Self-Organization*. Cambridge Univ. Press.

[24] SCHEIDEGGER, A. E. (1967). A stochastic model for drainage patterns into an intramontane trench. *Bull. Assoc. Sci. Hydrol.* **12** 15–20.

[25] SOUCALIUC, F., TÓTH, B. and WERNER, W. (2000). Reflection and coalescence between independent one-dimensional Brownian paths. *Ann. Inst. H. Poincaré Probab. Statist.* **36** 509–545. MR1785393

[26] TÓTH, B. and WERNER, W. (1998). The true self-repelling motion. *Probab. Theory Related Fields* **111** 375–452. MR1640799

[27] TROUTMAN, B. and KARLINGER, M. (1998). Spatial channel network models in hydrology. In *Stochastic Methods in Hydrology* (O. E. Barndorff-Nielsen, V. K. Gupta, V. Perez-Abreu and E. Waymire, eds.) 85–121. World Scientific, Singapore.

[28] TSIRELSON, B. (2004). Scaling limit, noise, stability. *Lectures on Probability Theory and Statistics*: *Ecole d'Eté de Probabilités de Saint-Flour XXXII. Lecture Notes in Math.* **1840** 1–106. Springer, Berlin. MR2079671

[29] VARADHAN, S. R. S. (1968). *Stochastic Processes*. Courant Institute of Math. Sciences, New York.

[30] WAYMIRE, E. (2002). Multiscale and multiplicative processes in fluid flows. In *Lectures on Multiscale and Multiplicative Processes in Fluid Flows. Centre for Mathematical Physics and Stochastics Lecture Notes* **11** 1–75. Dept. Math. Sciences, Univ. Aarhus.

[31] WINN, D. (2002). Gradient-directed diffusions and river network models. In *Lectures on Multiscale and Multiplicative Processes in Fluid Flows. Centre for Mathematical Physics and Stochastics Lecture Notes* **11** 149–154. Dept. Math. Sciences, Univ. Aarhus.

L. R. G. FONTES
INSTITUTO DE MATEMÁTICA E ESTATÍSTICA
UNIVERSIDADE DE SÃO PAULO
RUA DO MATÃO 1010
CEP 05508–090
SÃO PAULO SP
BRAZIL

M. ISOPI
DIPARTIMENTO DI MATEMATICA
UNIVERSITÀ DI ROMA "LA SAPIENZA"
P.LE ALDO MORO
00185 ROMA
ITALY

C. M. NEWMAN
COURANT INSTITUTE
  OF MATHEMATICAL SCIENCES
NEW YORK UNIVERSITY
NEW YORK, NEW YORK 10012
USA

K. RAVISHANKAR
DEPARTMENT OF MATHEMATICS
SUNY COLLEGE AT NEW PALTZ
NEW PALTZ, NEW YORK 12561
USA
E-MAIL: ravishak@newpaltz.edu